\documentclass{article}
\usepackage{amsmath,graphicx,mlspconf}
\usepackage{amssymb, amsthm}
\usepackage{optidef}
\usepackage{algorithm}
\usepackage{url}
\usepackage{comment}
\usepackage{hyperref}
\usepackage[usenames,dvipsnames]{color}
\usepackage{algorithmic}
\newtheorem{definition}{Definition}

\definecolor{brightpink}{rgb}{1.0, 0.0, 0.5}

\copyrightnotice{}

\copyrightnotice{Submitted to MLSP}

\copyrightnotice{}

\copyrightnotice{Submitted to MLSP}

\toappear{}

\title{Algorithms for Boolean matrix factorization\\ 
using integer programming}
\name{Christos Kolomvakis, Arnaud Vandaele, Nicolas Gillis\thanks{\noindent The authors acknowledge the support by the F.R.S.-FNRS and the FWO (EOS, O005318F-RG47). NG acknowledges the Francqui Foundation.}}
\address{Department of Mathematics and Operational Research \\ University of Mons, Rue de Houdain 9, 7000 Mons}

\begin{document}

\maketitle

\begin{abstract}
Boolean matrix factorization (BMF) approximates a given binary input matrix as the product of two smaller binary factors. 
As opposed to binary matrix factorization which uses standard arithmetic,  
BMF uses the Boolean OR and Boolean AND operations to perform matrix products, which leads to lower reconstruction errors.  
BMF is an NP-hard problem. In this paper, we first propose an alternating optimization (AO) strategy that solves the subproblem in one factor matrix in BMF using an integer program (IP). We also provide two ways to initialize the factors within AO. 
Then, we show how several solutions of BMF can be combined optimally using another IP. This allows us to come up with 
a new algorithm: it generates several solutions using AO and then combines them in an optimal way. 
Experiments show that our algorithms outperform the state of the art on medium-scale problems.  
\end{abstract}
\section{Introduction}
\label{sec:intro}

Low-rank matrix approximations (LRMAs) are popular methods in machine learning, and have successfully been  applied in a wide variety of applications such as document classification, community detection, hyperspectral unmixing and recommender systems, to cite a few; see, e.g.,~\cite{markovsky2012low, udell2016generalized, NMF_book}. 
LRMAs perform dimensionality reduction by approximating an input data matrix as the product  two factors of smaller sizes. 
Depending on the problem at hand, different matrix models can be considered. Examples include principal component analysis (PCA) and its variants such as sparse~\cite{zou2006sparse} and robust~\cite{candes2011robust} PCA, and nonnegative matrix factorization (NMF)~\cite{lee1999learning}. 
If the input matrix has elements in $\{0,1\}$, then it makes sense to impose the factors to have elements in $\{0,1\}$ as well, leading to binary matrix factorization (bMF) and Boolean matrix factorizations (BMF)~\cite{miettinen2008discrete, miettinen2008discrete, zhang2007binary}.  However, due to bMF using the standard addition and multiplication, the approximation will typically produce elements that are not in $\{0,1\}$. 
BMF uses the  Boolean OR and AND operations  to  guarantee the approximation to be binary, which improves the interpretability of the model. BMF is a difficult combinatorial problem, and many works have been dedicated to the computation of BMFs~\cite{Bayesian_BoolMF, miron2021boolean, dalleiger2022efficiently, ribeiro-PAKDD2016, generalized_miron}; see also~\cite{miettinen2020recent} for a recent survey. 

This paper proposes new ways to compute BMFs, and is organized as follows. 
In Section~\ref{sec:sec2}, we formally define BMF.  
In Section~\ref{sec:AO}, we describe our proposed alternating optimization (AO) algorithm for BMF where the subproblems are quadratic integer programs (IPs). We also provide two initialization strategies for the AO algorithm. 
In Section~\ref{sec:Comb}, we provide an IP formulation to optimally combine several BMF solutions. 
In Section~\ref{numexp}, we show that our proposed algorithms outperform the state of the art on four medium-scale real-world data sets.

\section{Boolean matrix factorization (BMF)}
\label{sec:sec2}

Let us first define the matrix Boolean product. 
\begin{definition}[Boolean product]
	Given two Boolean matrices, $\mathbf{W} \in \{0,1\}^{m \times r}$ and $\mathbf{H} \in \{0,1\}^{r \times n}$, their \textit{Boolean product} is denoted $\mathbf{W} \circ \mathbf{H} \in \{0,1\}^{m \times n}$ and is defined for all $i,j$ as 
	\begin{equation}\label{Bool_Prod}
		(\mathbf{W} \circ \mathbf{H})_{ij} = \bigvee_{k = 1}^r \mathbf{W}_{ik}\mathbf{H}_{kj},
	\end{equation}
	where $\vee$ is the logical OR operation 
 (that is, $0 \vee 0 = 0$, $1 \vee 0 = 1$, and $1 \vee 1 = 1$). 
 Interestingly, $\mathbf{W} \circ \mathbf{H} = \min(1,\mathbf{WH})$ where $\mathbf{WH}$ is the usual matrix product between $\mathbf{W}$ and $\mathbf{H}$. 
\end{definition}

We can now define the BMF problem. 
\begin{definition}[BMF] 
	Given a Boolean matrix $\mathbf{X} \in \{0,1\}^{m \times n}$ and a factorization rank $r$, BMF aims to find matrices $\mathbf{W} \in \{0,1\}^{m \times r}$ and $\mathbf{H} \in \{0,1\}^{r \times n}$ that solve 
	\begin{equation}
	\min_{\mathbf{W} \in \{0,1\}^{m \times r}, \mathbf{H} \in \{0,1\}^{r \times n}}	\|\mathbf{X} - \mathbf{W} \circ \mathbf{H} \|_F^2. 
		\label{BMF}
	\end{equation}	
\end{definition}
\noindent In \cite{miettinen2008discrete}, it is proven that not only solving \eqref{BMF}, 
but also approximating \eqref{BMF}, is NP-hard. 

BMF allows one to find subset of columns and rows of $\mathbf{X}$ that are highly correlated, since the entries equal to one in each binary rank-one factor $\mathbf{W}(:,k)\mathbf{H}(k,:)$ correspond to a rectangular submatrix of $\mathbf{X}$ that should contain many entries equal to one. 
Applications of BMF include role mining  \cite{Role_mining_BMF} and bioinformatics \cite{liang2020bem, haddad2018identifying}; see also the recent survey~\cite{miettinen2020recent}. 

\section{Alternating optimization (AO) for BMF} \label{sec:AO}

Most algorithms for LRMAs rely on iterative block coordinate descent methods: the subproblem in $\mathbf{H}$ is solved for $\mathbf{W}$ fixed, and vice versa. The reason is that these subproblems are typically convex. For BMF, this is of course not the case. However, the advances in IP solvers, such as Gurobi~\cite{gurobi}, allows one to tackle medium-scale problems efficiently.

\subsection{IP formulation for BMF subproblems}

Assuming $\mathbf{W}$ is fixed in~\eqref{BMF}, we would like to solve the following Boolean least squares (BoolLS) problem in $\mathbf{H}$, that is, solve 
\[
\min_{\mathbf{H} \in \{0,1\}^{r \times n}} 
\|\mathbf{X} - \min(1, \mathbf{WH}) \|_F^2. 
\] 
Because of the nonlinearity in the objective, this cannot be solved directly with standard IP solvers. 
Note that the problem in each column of $\mathbf{H}$ is independent: 
\begin{equation} \label{BoolLS} 
\min_{\mathbf{H}(:,j) \in \{0,1\}^{r \times n}} 
\| \mathbf{X}(:,j) - \min(1, \mathbf{WH}(:,j)) \|_F^2. 
\end{equation}
For simplicity, let $\mathbf{h} = \mathbf{H}(:,j)$ and \mbox{$\mathbf{x} = \mathbf{X}(:,j)$}. Given $\mathbf{W}$ and $\mathbf{x}$, we need to solve 
$\min_{\mathbf{h} \in \{0,1\}^r} 
\| \mathbf{x} - \min(1, \mathbf{Wh}) \|_F^2$. 
Introducing the variable $\mathbf{z} = \min(1, \mathbf{Wh})$, \eqref{BoolLS} can be reformulated as follows: 
\begin{equation} \label{bool_ls}
\min_{\mathbf{h} \in \{0,1\}^r, \mathbf{z} \in \{0,1\}^m}  \|\mathbf{x} - \mathbf{z}\|_F^2 
\quad 
\textrm{s.t.}  \quad 
\frac{\mathbf{Wh}}{r} \leq \mathbf{z} \leq \mathbf{Wh}. 
\end{equation} 
In fact, for $\mathbf{W}$, $\mathbf{h}$ and $\mathbf{z}$ binary, $\frac{\mathbf{Wh}}{r} \leq \mathbf{z} \leq \mathbf{Wh}$ if and only if $\mathbf{z} = \min(1,\mathbf{Wh})$, since  $\mathbf{Wh} \in \{0,1,\dots,r\}$. 
Now \eqref{bool_ls} is a convex quadratic optimization problem with linear constraints over binary variables. Such problems can be solved with commercial software, 
and  we make use of Gurobi~\cite{gurobi}. 

To give an idea of the computational time required for Gurobi to solve~\eqref{bool_ls}, let us perform the following experiment for various values of $m$ and $r$. The setting is as follows: we generate the entries of $\mathbf{W}$ and $\mathbf{h}$ using the uniform distribution in $[0,1]$, and then threshold all elements to convert them to binary. For all ranks tested, apart from $r = 50$, if an element is larger than 0.7, it is converted to 1, otherwise we convert it to 0. For $r = 50$, the threshold is set to 0.8. 
We chose relatively sparse $\mathbf{W}$ and $\mathbf{h}$ to make sure $\min(1,\mathbf{Wh})$ is not the all-one vector (in fact, we regenerate $\mathbf{W}$ and $\mathbf{h}$ if $\min(1,\mathbf{Wh})$ is the all-one or all-zero vector).  
Then we set $\mathbf{x} = \min(1,\mathbf{Wh})$, and 10\% of the entries of $\bf x$ are flipped randomly\footnote{The noiseless BoolLS problem is much easier to solve since $\mathbf{h}_k = 1$ if and only if the support of the $k$th column of $\mathbf{W}$ is contained in that of $\mathbf{x}$, that is, $\mathbf{W}(:,k) \leq \mathbf{x}$.}. Table~\ref{BoolLSnoisytime} reports the results. 
\begin{table}[H]
    \begin{center}
\begin{tabular}{||c |c c c c ||} 
 \hline
$r \backslash m$ & 100 & 1000 & 5000 & 10000\\ [0.5ex] 
 \hline\hline
 2 & 0.002 & 0.02 & 0.47 & 1.8\\ 
 \hline
 5 & 0.003 & 0.03 & 0.47 & 1.87 \\
 \hline
 10 & 0.005 & 0.04 & 0.56 & 2.2 \\
 \hline
 20 & 0.012 & 0.12 & 2.4 & 12.0 \\
 \hline
 50 & 0.049 & 2.52 & 38.0 & 235 \\  
 \hline
\end{tabular}
\caption{Average execution time in seconds over 30 trials of Gurobi to solve noisy BoolLS problems~\eqref{bool_ls} for various values of $m$ and $r$.}
 \label{BoolLSnoisytime}
\end{center}
\end{table} 
 The results are encouraging since solving~\eqref{bool_ls} can be done exactly with Gurobi in a reasonable amount of time, even for relatively large problems, e.g., it takes on average Gurobi 12 seconds to solve this problem with $m = 10^4$ and $r=20$. 
 Note that one could also use a timelimit for Gurobi, so that Gurobi would return the best solution found within the allotted time (often IP solvers take much more time to guarantee global optimality rather than finding the optimal solution).

\subsection{AO for BMF}

We can now solve BMF via AO over the factors $\mathbf{W}$ and $\mathbf{H}$ alternatively. Since $\|\mathbf{X} - \min(1,\mathbf{WH})\|_F^2
= \| \mathbf{X}^\top - \min(1,\mathbf{H}^\top \mathbf{W}^\top)\|_F^2$, the problem in $\mathbf{W}$ for $\mathbf{H}$ fixed has the same form. 
We update $\mathbf{H}$ in a column-by-column fashion by solving the independent BoolLS of the form~\eqref{BoolLS}, 
and similarly for $\mathbf{W}$ row by row. Algorithm~\ref{alg:bool_AO} summarizes the AO strategy. We have added a safety procedure within AO (steps~\ref{step5}-\ref{step9}): it may happen that some rows of $H$ are set to zero (for example, if $\mathbf W$ is not well initialized). 
In that case, we reinitialize these rows as the rows of the residual $\mathbf R = \max\big(0,\mathbf X - \max(1, \mathbf{W}_i \mathbf{H}_i)\big)$ whose entries have the largest sum. This guarantees that the error will decrease afer the update of $\mathbf{W}$.  


 



	

\begin{algorithm}[ht] 
\caption{AO algorithm for BMF - \texttt{AO-BMF}}
\begin{algorithmic}[1] \label{alg:bool_AO} 
\REQUIRE Input matrix $\mathbf{X} \in \{0,1\}^{m \times n}$, initial factor matrix $\mathbf{W}_0 \in \{0,1\}^{m \times r}$, maximum number of iterations maxiter. 

\ENSURE $\mathbf{W} \in \{0,1\}^{m \times r}$ and $\mathbf{H} \in \{0,1\}^{r \times n}$ such that $X \approx \min(1,WH)$. 
 
    \medskip  

\STATE $i = 1$, $e(0) = \|\mathbf{X}\|_F^2$, $e(1) = \|\mathbf{X}\|_F^2-1$. 

\WHILE{$e(i) < e(i-1)$ and $i \leq$ maxiter}
\STATE $\mathbf{H}_i = 
\texttt{BoolLS}(\mathbf{X},\mathbf{W}_{i-1})$. 
\STATE $\mathcal{K} = \{k \ | \ \mathbf{H}_i(k,:) = 0\}$. \label{step5}
\IF{ $\mathcal{K} \neq \emptyset$ }
    \STATE $\mathbf R = \max\big(0,\mathbf X - \max(1, \mathbf{W}_i \mathbf{H}_i)\big)$. 
    \STATE $\mathbf{H}_i(\mathcal{K},:) = \mathbf{R}(\mathcal{I},:)$, where $\mathcal{I}$ contains the indices 
    \STATE  of the $|\mathcal{K}|$ rows of $\mathbf{R}$ with largest sum. \label{step9} 
\ENDIF 
\STATE $\mathbf{W}_i = \texttt{BoolLS}(\mathbf{X}^\top, \mathbf{H}_{i}^\top)^\top$. 
\STATE $i = i + 1$. 
\STATE $e(i) = \| \mathbf{X} - \min(1, \mathbf{W}_i \mathbf{H}_i) \|_F^2$. 
\ENDWHILE
\RETURN $(\mathbf{W},\mathbf{H}) = (\mathbf{W}_i, \mathbf{H}_i)$. 
\end{algorithmic}  
\end{algorithm}

Typically, AO needs a very small number of iterations to converge, given the combinatorial nature of the problem. As we will report in Section~\ref{numexp}, among more than 6000 runs on 4 datasets with 3 different ranks, on average 3.7 iterations are needed, with a maximum of 14.


\subsection{Initialization of AO}

In this section, we provide two initialization strategies for AO-BMF, that is, Algorithm~\ref{alg:bool_AO}.

\paragraph*{Randomly selecting columns or rows of $\mathbf{X}$}

 AO-BMF only requires $\mathbf{W}$ to be initialized. By symmetry, it could also be initialized only with $\mathbf{H}$, starting the AO algorithm by optimizing over  $\mathbf{W}$. 
A simple, fast and meaningful strategy to initialize AO-BMF is to initialize $\mathbf{W}$ (resp.\ $\mathbf{H}$) with a subset of the columns (resp.\ rows) of $\mathbf{X}$, that is, 
set  
$\mathbf{W} = \mathbf{X}(:,\mathcal{K})$ (resp.\ $\mathbf{H} =  \mathbf{X}(\mathcal{K},:)$)  where $\mathcal{K}$ is a randomly selected set of $r$ indices of the columns (resp.\ rows) of $\mathbf{X}$. 

\paragraph*{NMF-based initialization} 

The second initialization we propose relies on NMF. NMF approximates $\mathbf{X}$ with $\mathbf{WH}$ where $\mathbf{W}$ and $\mathbf{H}$ are nonnegative. 
We use an NMF algorithm from~\url{https://gitlab.com/ngillis/nmfbook/} which itself initializes the entries of 
 $\mathbf{W}$ and $\mathbf{H}$ using the uniform distribution in $[0,1]$. Once an NMF solution is computed, we binarize it using the following two steps: 
 \begin{itemize}
     \item Normalize the columns of $\mathbf{W}$ and rows of $\mathbf{H}$ such that $\max(\mathbf{W}(:,k)) = \max(\mathbf{H}(k,:))$ for all $k$, using the scaling degree of freedom in NMF, that is,  $\mathbf{W}(:,k) \mathbf{H}(k,:) = (\alpha \mathbf{W}(:,k))(\alpha^{-1} \mathbf{H}(k,:))$ for $\alpha > 0$. 
     
     \item Set the entries of $\mathbf{W}$ and $\mathbf{H}$ to 0 or 1 using a given threshold $\delta$ which is generated uniformly at random  in the interval $[0.3, 0.7]$.  An alternative would be to use a grid search approach to determine $\delta$, similar to~\cite{zhang2007binary, truong2021boolean}. 
However, we observed that selecting $\delta$ randomly performs better on average.  
 \end{itemize}

\section{Combining multiple BMF solutions}  \label{sec:Comb}

Algorithm~\ref{alg:bool_AO}, AO-BMF, is able to relatively quickly generate locally optimal solutions for~(\ref{BMF}), in the sense that they cannot be improved by optimizing $\textbf{W}$ or $\textbf{H}$ alone. 
A first natural approach to generate good solutions to BMF is using  
multiple initializations, and keeping the best solution. We will refer to this strategy as multistart AO (MS-AO).  

However, it is possible to combine a set of solutions in a more effective way. 
Assume we have generated $p$ rank-$r$ BMFs: 
$\mathbf{W}_1 \mathbf{H}_1, \dots, \mathbf{W}_p \mathbf{H}_p$. This gives $rp$ binary rank-one factors, namely $\mathbf{W}_\ell(:,k)\mathbf{H}_\ell(k,:)$ for $k=1,\dots,r$ and $\ell=1,\dots,p$. Let us denote these rank-one binary matrices $\mathbf{R}_i$ for $i=1,\dots,N$ with\footnote{In practice, we delete duplicated rank-one factors so that $N \ll rp$.} $N = rp$. 
To generate a better rank-$r$ BMF, we can pick $r$ rank-one binary factors among the $\mathbf{R}_i$'s, by solving the following combinatorial problem: 
\[
\min_{\mathbf{y} \in \{0,1\}^N} \Big\|\mathbf{X} - \min\big( 1 , \sum_i \mathbf{y}_i \mathbf{R}_i \big) \Big\|_F^2 \; \text{ such that } \; \sum_i y_i = r.  
\] 
The variable $y \in \{0,1\}^N$ encodes the $r$ selected rank-one factors, that is, $y_i = 1$ if $\mathbf{R}_i$ is selected in the BMF. 
As for BoolLS, we can reformulate this problem as a quadratic IP: 
        \begin{align} \label{Comb_BoolMF}
 \min_{\mathbf{y} \in \{0,1\}^N, 
 \mathbf{Z} \in \{0,1\}^{m\times n}} 
            &  \left\| \mathbf{X} - \mathbf{Z} \right\|_F^2 \\
 \text{ such that } \quad 
 \sum_{i=1}^N y_i = r, 
 & \frac{\sum_{i=1}^N y_i \mathbf{R}_i}{r} \leq \mathbf{Z} \leq \sum_{i=1}^N y_i \mathbf{R}_i. \nonumber 
        \end{align} 
        Denoting $y_i^*$ the optimal solution of~\eqref{Comb_BoolMF},  
the rank-$r$ BMF obtained, $\sum_i y_i^* \mathbf{R}_i$, is guaranteed to be at least as good as all the solutions $\{\mathbf{W}_i \mathbf{H}_i\}_{i=1}^p$, since they are feasible solutions of~\eqref{Comb_BoolMF}. 

Gurobi can solve medium-scale problem of the form~\eqref{Comb_BoolMF} in a reasonable amount of time. Table~\ref{runtimecomb} reports the time to solve~\eqref{Comb_BoolMF} for $m=n=101$, $r=5$, and $N=2^k \ 50$ for $k=0,1,\dots,4$. For example, to combine 160 rank-5 solutions (hence 800 rank-one factors) of a 101-by-101 matrix, it takes about 30 seconds. 
\begin{table}[H]
	\begin{center}
		\begin{tabular}{||c | c c c c c||} 
			\hline
    $N$         &   50 & 100    & 200 & 400  & 800\\ \hline
 time (s.)    &  4.0 & 4.3  & 7.6 & 14.2 & 30.6 \\ \hline
		\end{tabular}
	\end{center}
	\caption{Run times to combine $N$ rank-one factors for a 101-by-101 matrix with $r=5$ (namely, the apb data set, see Section~\ref{numexp}).}
	\label{runtimecomb}
\end{table}

In practice, if $N$ is too large, we do not need to take into account all rank-one factors, and can only consider rank-one factors corresponding to the best BMFs. 

Finally, once a solution combining several BMFs is computed, we will further improve it using AO. 

We will refer to this algorithm, namely generating $p$ solutions with AO-BMF, then combining them solving~\eqref{Comb_BoolMF}, and then applying AO to that solution as MS-Comb-AO.





\section{Numerical Experiments}  \label{numexp}

All experiments are performed with a 
12th Gen Intel(R) Core(TM) i9-12900H  2.50 GHz, 32GB RAM, 
on MATLAB R2019b. 
The code and data sets are available on \url{https://gitlab.com/ngillis/BooleanMF} 

\paragraph*{Data sets} We will perform experiments on four real data sets used in~\cite{BoolMF_IP}, and which come from~\cite{UCI, krebs2008network}; see Table~\ref{datasets}. \vspace{-0.5cm}
\begin{table}[H]
	\begin{center}
		\begin{tabular}{||c | c c c c||} 
			\hline
			& zoo & heart & lymp & apb \\ [0.5ex] 
			\hline
			$m \times n$ & $101 \times 17$ & $242 \times 22$ & $148 \times 44$ & $105 \times 105$\\ 
			\hline
		\end{tabular}
	\end{center}
	\caption{Four binary real-world data sets.}
	\label{datasets}
\end{table}
\noindent As in \cite{BoolMF_IP}, we use $r= 2, 5, 10$ for all data sets. 
In~\cite{BoolMF_IP}, authors proposed two non-trivial IP-based approaches for BMF that perform well against the state of the art (they used a 20 minutes time limit for their method), 
namely against 
a greedy scheme~\cite{BoolMF_IP}, 
ASSO and ASSO++~\cite{miettinen2008discrete}, 
a penalty formulation from~\cite{zhang2007binary}, 
and an NMF-based heuristic. 
Table~\ref{bestres} reports the best result of all these methods for the four data sets. From now on, we will report the result of a BMF by comparing it with the best solution in Table~\ref{bestres}. 
\begin{table}[H]
	\begin{center}
		\begin{tabular}{||c | c |c |c||} 
			\hline
			& $r = 2$ & $r = 5$ & $r = 10$ \\ [0.5ex] 
			\hline
			zoo &  271 & 	126	 & 39\\
			\hline 
			heart &  1185 &	737 &	419 \\
			\hline
			lymp &  1184	& 982 &	728 \\ 
			\hline
			apb &  776 &	684	& 573\\
			\hline
		\end{tabular}
	\end{center}
	\caption{Objective function 
 $\|\mathbf{X}-\min(1,\mathbf{WH})\|_F^2$ of the best solution found by various algorithms in~\cite[Table~4, page~20]{BoolMF_IP}.} 
	\label{bestres}
\end{table}

\paragraph*{Results from two other papers} We have also run two recent algorithms, from \cite{dalleiger2022efficiently} and \cite{generalized_miron}, on these data sets. 
The method in~\cite{dalleiger2022efficiently}  is based on a continuous reformulation of BMF and uses alternating proximal projected gradient method to solve it. We have used the parameters of the algorithm recommended by the authors. 
Table~\ref{nips_res} reports the results. 
{\color{black} Out of curiosity, we initialized AO with the best solutions found. We observe in Table~\ref{nips_res} that AO is able to significantly improve these solutions, showing that the continuous reformulations of \cite{dalleiger2022efficiently} is not able to generate locally optimal solutions.} 

\begin{table}[H]
	\begin{center}
		\begin{tabular}{||c | c |c |c||} 
			\hline
			& $r = 2$ & $r = 5$ & $r = 10$ \\ [0.5ex] 
			\hline
			zoo &  +11 $ | $ \ \ 0 \ \  & +5 \ $ | $ \ -1 & +18 \ $ | $ \ +5\\
			\hline 
			heart & +65 $ | $ +18  & +29 $ | $ \ -1 \ & \ +33 \ $ | $ +26\\
			\hline
			lymp &  +64 $ | $ +25 & +88 $ | $ -10 & +203 $ | $ +60\\ 
			\hline
			apb &  +72 $ | $ +30 & \ +59 $ | $ +43  & \ +33 \ $ | $ +15\\
			\hline
		\end{tabular}
	\end{center}
	\caption{Best result obtained with the method in \cite{dalleiger2022efficiently} (left), and result of this best solution improved by AO (right).
  The numbers indicate the difference compared with the best values in Table~\ref{bestres}. 
  A negative value means an improvement, a positive value means a worse solution. 
  } 
	\label{nips_res}
\end{table}

Next, we show the results for the method in \cite{generalized_miron}. Note that this paper provides a rather general approach, allowing for other Boolean operations between the factors to approximate the entries of $\mathbf{X}$. Their approach is also based on some continuous reformulations, and the use of gradient and descent methods. Table~\ref{miron_res} reports their result, using 15 trials where each trial uses 10 random initialization, each optimized with 2000 iterations. {\color{black} We also use AO to improve the best solutions found by this method, and observe a similar behavior as for the method from~\cite{dalleiger2022efficiently}.}
\begin{table}[H]
	\begin{center}
		\begin{tabular}{||c | c |c |c||} 
			\hline
			& $r = 2$ & $r = 5$ & $r = 10$ \\ [0.5ex] 
			\hline
			zoo &  \ +2 \ \ $ | $ \ 0 & \ +7 \ $ | $ \ +4 & +10 \ $ | $ \ +4\\
			\hline 
			heart & +165 $ | $ +9 & \ +33 $ | $ \ +4 \ & \ +86 \ $ | $ +57\\
			\hline
			lymp & +214 $ | $ \ -3 & +49 $ | $ -16 & +123 $ | $ \ +5\ \  \\ 
			\hline
			apb & +48 \ $ | $ \ 0 & \ +51 $ | $ +23& +112 $ | $ +43\\
			\hline
		\end{tabular}
	\end{center}
	\caption{Best result obtained with the method in \cite{generalized_miron} (left), and result of this best solution improved by AO (right).  
  The numbers indicate the difference compared with the best values in Table~\ref{bestres}.
 }
	\label{miron_res}
\end{table}

In summary, on these data sets, we observe that the algorithms from~\cite{dalleiger2022efficiently, generalized_miron} provide solutions which are much worse than the best solutions provided in~\cite{BoolMF_IP}, and that AO can considerably improve the solutions of these two methods, sometimes obtaining better solution than the best from~\cite{BoolMF_IP} (e.g., for the lymp data set with $r=5$).  

\paragraph*{Results of our proposed algorithms} 

Let us provide the results for our two algorithms. 

\noindent \textbf{[1.] Multiple starts of AO-BMF (MS-AO).} We generate as many BMFs as possible with AO-BMF within the time $T$, and return the best solution. The initialization of AO-BMF is chosen alternatively as one of the two strategies (NMF-based or random columns/rows of $\bf{X}$).  
Table~\ref{aoresults} reports the results for $T$ equal to 30 seconds and 5 minutes. 
\begin{table}[H]
	\begin{center}
		\begin{tabular}{||c | c |c |c||} 
			\hline
			& $r = 2$ & $r = 5$ & $r = 10$ \\ [0.5ex] 
			\hline
zoo & \ 0 \ $ \ | \ $ 0  &  -1  \ $ \ | \ $ -1   & \ +3 \ $ \ | \ $ \ \ 0  \ \ \\ \hline 
heart & \ +2 \ $ \ | \ $ +2   &  -1 \  $ \ | \ $ -1  
   & 0 \ \ $ \ | \ $  \ \ 0    \\ 
			\hline
lymp & \ \ \ -10  \  $ \ | \ $ -10 \ \  &  -25 \ $ \ | \ $ -32  & -15 $ \ | \ $ -34 \\
			\hline
apb & \ 0 \ $ \ | \ $ 0 & +6 \ $ \ | \ $ \ -6  & +4 $ \ | \ $ +2  \\ 
			\hline
		\end{tabular}
	\end{center}
	\caption{Results of the MS-AO strategy for 30 seconds (left) and 5 minutes (right).  The numbers indicate the difference compared with the best values in Table~\ref{bestres}.}   
	\label{aoresults}
\end{table} 
Quite surprisingly, MS-AO is already able to perform on par or improve upon the state of the art. With a 30 seconds time limit, it does on 8 out of 12 cases: 5/12 cases with improvements, sometimes significant as for the lymp data set, and 3/12 cases with the same objective. 
However, for 4 data sets, it is not able to achieve the best solution reported in Table~\ref{bestres}, although the generated solutions are only slightly worse (+6 at most).  
With a 5 minutes time limit, it performs better or on par on 10 out of the 12 cases, and it produces a slightly worse solution in two cases (+2). 

To give an idea of the performance of Gurobi on these medium-scale problems, Table~\ref{iterBoolMFgurobi} reports the number of BMFs generated within 5 minutes for each data set, as well as the average number of iterations required for AO to converge. 
The average number of iterations for AO to converge is 3.7, and the largest number of iterations AO needed among these 6165 BMFs is 14.   
\begin{table}[H]
	\begin{center}
		\begin{tabular}{||c | c |c |c||} 
			\hline
			& $r = 2$ & $r = 5$ & $r = 10$ \\ [0.5ex] 
			\hline
			zoo & 844 $ \ | \ $ 3.3  & 703 $ \ | \ $ 3.2   & 596 $ \ | \ $ 3.3\\   			\hline  
			heart &  503   $ \ | \ $ 3.0   & 280     $ \ | \ $ 3.5  &  484     $ \ | \ $ 3.5 \\
			\hline
			lymp & 558 $ \ | \ $ 3.3  &   568     $ \ | \ $ 4.2  &  269     $ \ | \ $ 4.8  \\
			\hline
			apb & 591 $ \ | \ $ 3.7  & 381    $ \ | \ $ 4.2  &  388   \  $ \ | \ $ \ 4.8 \\ 
			\hline
		\end{tabular}
	\end{center}
	\caption{Number of BMFs generated via AO within 5 minutes  (left), and average number of iterations required for AO to converge (right).}
	\label{iterBoolMFgurobi}
\end{table}

\noindent \textbf{[2.] Multiple starts, combination and AO (MS-Comb-AO).} As explained in Section~\ref{sec:Comb}, we generate as many BMFs as possible with AO-BMF (as for AO-MS) within time $3T/4$, 
and then combine them by solving~\eqref{Comb_BoolMF} with a time limit of $T/4$. We used the same random seed so that the solutions generated are the same as for  AO-MS, except that less solutions are generated, since only 3/4 of the total time is spent for that.  
Table~\ref{result_comb} reports the results. 
\begin{table}[H]
	\begin{center}
		\begin{tabular}{||c | c |c |c||} 
			\hline
& $r = 2$ & $r = 5$ & $r = 10$ \\ [0.5ex] 
			\hline
zoo & \ 0 \  $ \ | \ $ 0   & -1 $ \ | \ $ -1 & 0  $ \ | \ $ 0 \\ 
			\hline 
heart &  +2  $ \ | \ $ 0  & -1 $ \ | \ $ -1  & 0 $ \ | \ $ 0  \\
			\hline
lymp &  \ \ -10   $ \ | \ $ -10  & -25 $ \ | \ $ -32 & -15 $ \ | \ $ -34 \\
			\hline
apb & \ 0 \ $ \ | \ $ 0  & -1 $ \ | \ $ -6   & +2 $ \ | \ $ -7  \\ 
			\hline
		\end{tabular}
	\end{center}
	\caption{Results of the MS-Comb-AO algorithm for 30 seconds (left) and 5 minutes (right).  The numbers indicate the difference compared with the best values in Table~\ref{bestres}.} 
	\label{result_comb}
\end{table}

For some data set, AO is already able to generate very good solutions, and hence solving~\eqref{Comb_BoolMF} is not useful, e.g., for the lymp data set.  
However, for most cases, this combination is beneficial, sometimes significantly. In particular, with the timelimit of 30 seconds on the apb dataset for $r=5$, the best solution found by AO-MS  has error 689 
(+6) while the combination leads to an error of 682 
(-1); for the zoo data set with $r=10$, 
it goes from 42 (+3) to 39 (0). A similar behavior is observed for 5 minutes: for the apb data set with $r=10$  
from +2 to -7, and the heart data set from +2 to 0. 

To conclude,  MS-Comb-AO with a 5-minute timelimit is able to either perform on par with the state of the art (we suspect that for these data sets, the corresponding solutions are optimal), or outperform it, sometimes significantly (in particular for the lymp data sets).

\begin{figure}[ht!]
	\begin{center} 
		\includegraphics[width=0.4\textwidth]{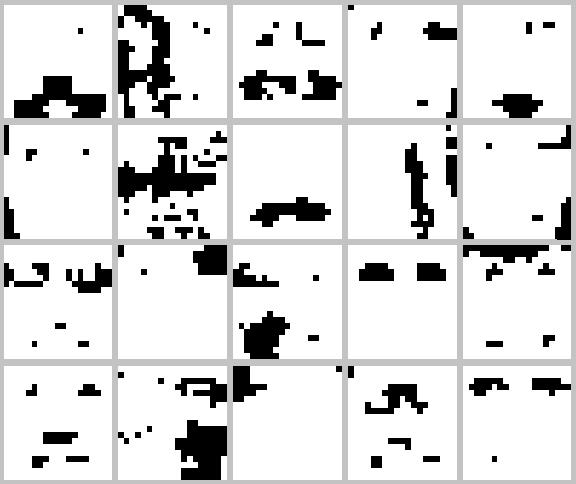}
		\caption{Facial features extracted by AO-BMF on the CBCL data set.}  
		\label{cbclim}
	\end{center}
\end{figure}

\paragraph*{Experiment on facial images} As a last experiment, let us apply AO-BMF to a larger data set, the well-known CBCL facial images. It was used in the seminal paper of Lee and Seung to show that NMF is able to extract facial features~\cite{lee1999learning}. 
Let us apply  AO-BMF on this very same data set. Each column of the data matrix $\mathbf{X} \in \mathbb{R}^{361 \times 2429}$ contains a vectorized facial image of size $19 \times 19$, and is not binary but satisfies $\mathbf{X}(i,j) \in [0,1]$ for all $(i,j)$. Quite interestingly, AO-BMF can be applied to any input matrix $\mathbf{X}$, even if it is not binary. We use the NMF-based initialization for AO-BMF with $r=20$, and it converges in 15 iterations within about 1 hour. AO-BMF provides a binary approximation of $\mathbf{X}$ with relative error $\frac{\|\mathbf{X} - \min(1,\mathbf{WH})\|_F}{\|\mathbf{X}\|_F} = 55.94\%$. 
Figure~\ref{cbclim} displays the meaningful binary facial features extracted by AO-BMF as the columns of $W$. 

\section{Conclusion and further work} 

In this paper, we have designed an alternating optimization (AO) strategy to tackle Boolean matrix factorization (BMF) using interger programming (IP). We have also shown how to combine several solutions in an optimal way, using IP as well. We finally showed how these two strategies are able to outperform the state of the art on 4 real-world medium-scale data sets. 

Further work includes the adaptation of our algorithms for large-scale data sets. A simple, yet possibly effective approach, is to use a time limit for Gurobi to tackle the IP subproblems. 
Typically, Gurobi is able to produce quickly high quality solutions. 
The design of fast heuristic algorithms to solve the Boolean least squares problem~\eqref{bool_ls} would also be useful, since we have validated the effectiveness of AO for BMF. 
Another line of work, that we are currently exploring, is to design more sophisticated combination strategies.  
This has allowed us  to generate even better solutions for three cases: 

 \noindent \quad $\bullet$ for lymp with $r=5$, a BMF with error 939 (-43), 
 
 \noindent \quad $\bullet$ for lymp with $r=10$, a BMF with error 680 (-48). 
 
 \noindent \quad $\bullet$ for apb with $r=5$, a BMF with error 677 (-7).

\noindent We will present these results in an extended version of this paper.


\paragraph*{Acknowledgments} We thank Sebastian Miron for sending us his BMF code~\cite{generalized_miron}, and Sebastian Dalleiger for providing us with good parameters for his BMF algorithm~\cite{dalleiger2022efficiently}.

\bibliographystyle{IEEEbib}
\bibliography{References}

\end{document}